\documentclass[submission,copyright,creativecommons]{eptcs}
\providecommand{\event}{GASCOM 2024} 

\usepackage{iftex}

\ifpdf
  \usepackage{underscore}         
  \usepackage[T1]{fontenc}        
\else
  \usepackage{breakurl}           
\fi

\usepackage{amsmath, amsfonts, amsthm, amssymb}
\usepackage{latexsym, mathrsfs} 
\usepackage{mathtools}
\usepackage{enumitem}

\usepackage{tikz}
\usepackage{tikz-cd}
\tikzcdset{scale cd/.style={every label/.append style={scale=#1},
    cells={nodes={scale=#1}}}}
\usepackage{ytableau}
\usepackage[edges]{forest}
\usepackage[height=\baselineskip]{mahjong}

\usepackage{genyoungtabtikz}
\newcommand\Ygr{\Yfillcolour{gray!20}}
\newcommand\Yw{\Yfillcolour{white}}

\usepackage[capitalise]{cleveref}

\newcommand{\N}{\mathbb{N}}

\newcommand{\A}{\mathcal{A}}
\newcommand{\B}{\mathcal{B}}

\newcommand{\R}{\mathcal{R}}

\newcommand{\J}{\mathscr{J}}

\newcommand{\M}{\mathbf{M}}

\newcommand{\chinese}{\mathbf{Ch}}
\newcommand{\plax}{\mathbf{Plax}}

\newcommand{\Com}{\mathbf{Com}}
\newcommand{\Sylv}{\mathbf{Sylv}}

\renewcommand{\to}{\rightarrow}

\newcommand{\cont}{\operatorname{cont}}
\newcommand{\ev}{\operatorname{ev}}
\newcommand{\RW}{\operatorname{R}}
\newcommand{\dom}{\operatorname{dom}}
\newcommand{\dyck}{\mathbf{Dyck}}

\newcommand{\f}{\raisebox{-4px}{\mahjong{2p}}}

\newcommand{\emptyf}{\raisebox{-4px}{\mahjong{5z}}}

\newtheorem{thm}{Theorem}[section]
\newtheorem{prop}[thm]{Proposition}

\newtheorem{defi}[thm]{Definition}
\numberwithin{equation}{section}

\title{Power Quotients of Plactic-like Monoids}
\author{Antoine Abram\thanks{Antoine Abram was partially supported by CRSNG BESC D grant,}
\institute{LACIM, Université du Québec à Montréal,\\ Montréal, Québec}
\email{\href{mailto:abram.antoine@courrier.uqam.ca}{abram.antoine@courrier.uqam.ca}}
\and
Florent Hivert\thanks{Florent Hivert was partially supported by OpenDreamKit Horizon 2020 ERI (\#676541)}
\institute{LISN, Université Paris-Saclay,\\ Orsay, France}
\email{\href{mailto:florent.hivert@lri.fr}{florent.hivert@lri.fr}}
\and
James D. Mitchell\thanks{James~D.~Mitchell was partially supported by OpenDreamKit Horizon 2020 ERI (\#676541)}
\institute{School of Mathematics and Statistics, University of St Andrews,\\ St Andrews, Scotland}
\email{\href{mailto:jdm3@st-andrews.ac.uk}{jdm3@st-andrews.ac.uk}}
\and
Jean-Christophe Novelli\thanks{Jean-Christophe~Novelli was partially supported by the ANR program CARPLO.}
\institute{LIGM, Université Gustave-Eiffel,\\ Marne-la-Vallée, France}
\email{\href{mailto:novelli@univ-mlv.fr}{novelli@univ-mlv.fr}}
\and
Maria Tsalakou,
\institute{School of Mathematics and Statistics, University of St Andrews,\\ St Andrews, Scotland}
\email{\href{mailto:mt200@st-andrews.ac.uk}{mt200@st-andrews.ac.uk}}
}
\def\titlerunning{Power Quotients of Plactic-like Monoids}
\def\authorrunning{A.~Abram, F.~Hivert, J.~D.~Mitchell, J.-C.~Novelli, M.~Tsalakou}

\begin{document}
\maketitle

\begin{abstract}
  In this paper we describe the quotients of several plactic-like monoids by the
  least congruences containing the relations $a ^ {\sigma(a)} = a$ with
  $\sigma(a)\geq 2$ for every generator $a$.
  The starting point for this description is the recent paper of
  Abram and Reutenauer about the so-called \emph{stylic monoid} which happens to
  be the quotient of the plactic monoid by the relations $a ^ 2 = a$ for every
  letter $a$. The plactic-like monoids considered are the plactic monoid
  itself, the Chinese monoid, and the sylvester monoid.
  In each case we describe: a set of normal forms, and the idempotents; and obtain
  formulae for their size.
\end{abstract}

\section{Introduction}
\label{section-introduction}

The plactic monoid $\plax(\A)$ over an ordered alphabet $(\A,<)$, can be defined
as the quotient of the free monoid $\A^*$ on $\A$
by
identifying words that produce the same Young tableau using Robinson-Schensted
insertion algorithm~\cite{R38, S61}.

Knuth~\cite{K70} found
an explicit presentation as the quotient of $\A^*$ by the relations:
\[
  acb  =  cab \quad \text{if } a \leq b < c \quad\text{and}\quad
  bac =  bca \quad \text{if } a < b \leq c, \quad \text{with }a, b, c\in \A. 
\]
For more details see \cite{LS81}, or \cite[Chapter 5]{L02}.
Due to its link with symmetric functions and representation
theory, the plactic monoid is a central object in algebraic combinatorics
that has been widely studied in the literature.

Other monoids, whose relations are delineated in terms of insertion algorithms
on certain combinatorial objects, are often referred to as "plactic-like" monoids.
They exhibit a rich combinatorial structure and have applications in several topics including
geometry and representation theory.

Among others, this family contains the Chinese monoid $\chinese(\A)$~\cite{CEKNH01},
that has applications on Hecke atoms and the Bruhat order
(see \cite{HMP17}),
and the sylvester monoid $\Sylv(\A)$~\cite{HNT05},
which is related to the associahedra and the Loday-Ronco algebra of trees.

Recently, Reutenauer and the first author~\cite{AR22} discovered that the quotient of the
plactic monoid by the relations $a^2=a$, for every letter $a$,
has several interesting properties.
Inspired by their results,
we investigate more general finite quotients of plactic-like monoids.
For such monoid
$\M(\A)$ defined over an alphabet $\A$ and a function $\sigma:\A\to \N_{\geq
2}$, we study the quotients $\M(\A, \sigma)$ by the relations $a^{\sigma(a)} =
a$ for every $a\in \A$.
It turns out that the monoids are of two different
non-disjoint types. For the first type, that includes the plactic, Chinese, and
hypoplactic monoids,
$\M(\A, \sigma)$
can be naturally embedded in the cartesian product of $\M(\A, 2)$ with the
commutative monoid $\Com(\A, \sigma)$.
The second type, that includes the
hypoplactic and the sylvester monoids, have words with a particular property in
every equivalence class, that provides a set of normal forms.
In both types,
$\M(\A, 2)$ plays an important role in the structure of $\M(\A,\sigma)$
for any $\sigma$.
In addition, $\M(\A, 2)$ has a rich combinatorial structure usually related to the one of $\M(\A)$.

For the first type, our knowledge of the stylic monoid helps us
in understanding the structure of $\plax(\A, \sigma)$.
We then consider another example, namely the Chinese monoid,
by studying its $2$-quotient, which involves rich combinatorial objects,
and transpose this to the study of its general $\sigma$-quotient.

For the second type, we focus on the sylvester monoid and more particularly on
its 2-quotient.

\section{Words, Monoids and $\sigma$-Quotients}
\label{section-notation}

Let $\A$ be a finite ordered alphabet, $\A^*$ the free monoid over $\A$,
and $\Com(\A)$ the free commutative monoid over $\A$.

For $W\in \A ^*$, the \emph{content} of $W$ is the set of distinct
letters occurring in $W$, and is denoted by $\cont(W)$.
We call the natural surjection $\ev : \A^* \to \Com(\A)$ the \emph{evaluation map}
and say $\ev(W)$ is the \emph{evaluation} of the word $W$, for all $W \in \A^*$.
For a word $W= w_1 \cdots
w_n\in \A^ *$, an \emph{inflation} of $W$ is any
word of the form $w_1 ^ {\epsilon_1} \cdots w_n ^ {\epsilon_n}\in \A ^ *$ for
some $\epsilon_1, \ldots, \epsilon_n\in \N_{\geq 1}$.

Let $\M(\A)$ be a monoid defined by the presentation $\langle \A |
\R\rangle$ for some alphabet $\A$ and some relations
$\R\subseteq \A ^ * \times \A ^ *$. If $U, V\in \A^*$ have the same
image under the natural surjective homomorphism $\pi_{\M}: \A ^ * \to \M$,
we say that $U$ and $V$ are \emph{equivalent} in $\M(\A)$, and we write $U
\equiv_{\M(\A)} V$.

If $\R\subseteq \A ^* \times \A ^ *$ has the property that $\ev(U) =
\ev(V)$ for all $(U, V) \in \R$, then we say that the presentation
$\langle \A | \R\rangle$ is \emph{evaluation-preserving},
and, by extension, that $M(\A)$ is an \emph{evaluation-preserving monoid}.

\begin{defi}\label{def:Quotient}
Let $\M(\A)$ be an evaluation preserving monoid and $\sigma : \A \to \N_{\geq 2}$.
We define $\M(\A, \sigma)$ to be the
quotient of $\M(\A)$ obtained by adding the extra relations $(a ^ {\sigma(a)}, a)$
for every $a\in \A$ to the presentation $\langle \A | \R\rangle$.
If $\sigma:\A \to \N_{\geq 2}$ is constant with value $n$, then we write
$\M(\A,n)$ instead of $\M(\A, \sigma)$.
\end{defi}

Two types of monoids arise from the study of these quotients. These two types
are not mutually exclusive; the hypoplactic monoid is
of both types.

\section{Monoids of Type 1}
\label{s:type1}

Let $M(\A)$ be an evaluation-preserving monoid. For any $\sigma:\A\to\N_{\geq 2}$,
let $\theta:\M(\A,\sigma)\to\M(\A,2)$ and $\ev_\sigma:\M(\A,\sigma)\to\Com(A,\sigma)$
be the natural surjective morphisms.
We let $\phi_\sigma$ be the product map $\theta \times ev_\sigma$.

\begin{defi}\label{def:type1}
An evaluation-preserving monoid $\M(\A)$ is of \emph{type 1} if
for any $\sigma:\A\to \N_{\geq2}$, $\phi_\sigma$ is an embedding.
\end{defi}

\subsection{The $\sigma$-Plactic Monoids}
\label{section-plax}

In \cite{AR22}, the monoid which we define here as $\plax(\A, 2)$ was introduced
as the stylic monoid. We recall some necessary definitions and
refer the reader to \cite{AR22} for more details.

An \emph{$N$-tableau} is a semi-standard tableau such that its rows are
strictly increasing; and each row is contained in the row underneath. The
\emph{row reading} of a tableau $T$ is the word $\operatorname{R}(T)$ obtained
by reading each row from left to right and top to bottom. The $N$-tableaux have
an insertion algorithm, denoted $T \leftarrow W$, called the \emph{$N$-insertion}
that is similar to the Robinson-Schensted insertion in Young tableaux.

Using these $N$-tableaux and their insertion algorithm, one can prove
that the plactic monoid is of type 1.

Given $(T, e)\in \plax(\A, \sigma)\times \Com(\A, \sigma)$ such that
$\cont(T) = \cont(e)$, we define the \emph{row reading} $\RW_{\sigma}(T, e)
\in \A^*$ via the following algorithm starting with $\RW(T)$. For every $a\in
\cont(T)$, if $\alpha$ is the number of occurrences of $a$ in $\RW(T)$ and
$\beta$ is the exponent of $a$ in $e$, then we replace the last occurrence of
$a$ in $\RW(T)$ by $a^\gamma$ where $\gamma\in \N_{\geq 0}$ is the least
value such that $\gamma + \alpha = \beta \pmod{\sigma(a)-1}$.
For example, if $\sigma(x)=4$ for all $x\in\A$ then
\[
\RW_\sigma\left(\vcenter{\hbox{\gyoung(b,ab)}}, a^1b^1\right) = bab ^ 3
\qquad\text{and}\qquad
\RW_\sigma\left(\vcenter{\hbox{\gyoung(bc,abc)}}, a^3b^2c^1\right) = bca^3bc^3\,.
\]
Such a row reading $\RW_{\sigma}(T, e)$ is a preimage of $(T, e)$ under $\phi$;
and these row readings constitute a set of normal forms for
$\plax(\A, \sigma)$.

\subsection{The $\sigma$-Chinese Monoids}
\label{section-2-chinese}

The \emph{Chinese monoid}
$\chinese(\A)$ is defined by the presentation with generating set $\A$ and
relations \cite{CEKNH01}: for $a, b, c\in \A$,
\[
    cba =   cab = bca \quad\text{if } a < b < c,  \quad
    aba  =  baa, \quad bba = bab \quad \text{if } a < b.
\]

\begin{figure}[htb]
\centering  
  \scalebox{0.8}{
  $
\begin{array}{ccc}
\Yboxdim{1cm}
\begin{tikzpicture}[scale=0.5,baseline={(current bounding box.east)}]
  \tgyoung(0cm,0cm,%
  ::::!\Ygr1!\Yw:a,%
  :::!\Ygr;2;!\Yw:b,%
  ::!\Ygr3;;2!\Yw:d,%
  :!\Ygr;;3!\Yw;;:f,%
  !\Ygr2;1!\Yw;;;:g,%
  :g:f:d:b:a)
  \draw[ultra thick](0,-4)--(2,-4)--(2,-3)--(3,-3)--(3,-2)--(5,-2)--(5,1);
\end{tikzpicture}
&
\Yboxdim{1cm}
\begin{tikzpicture}[scale=0.5,baseline={(current bounding box.east)}]
  \tgyoung(0cm,0cm,%
  ::::!\Ygr;!\Yw:a,%
  :::!\Ygr;1;!\Yw:b,%
  ::!\Ygr;;;1!\Yw:d,%
  :!\Ygr;;1!\Yw;;:f,%
  !\Ygr;;1!\Yw;;;:g,%
  :g:f:d:b:a)
  \draw[ultra thick](0,-4)--(2,-4)--(2,-3)--(3,-3)--(3,-2)--(5,-2)--(5,1);
\end{tikzpicture}
      &
        \qquad
\begin{pmatrix}
  a & b & d & f & g \\
  a & a & a & d & f \\
\end{pmatrix}
\quad
\Yboxdim{1cm}
\begin{tikzpicture}[scale=0.5,baseline={(current bounding box.east)}]
  \tgyoung(0cm,0cm,%
  ::::!\Ygr;1!\Yw:a,%
  :::!\Ygr;;1!\Yw:b,%
  ::!\Ygr;;;1!\Yw:d,%
  :!\Ygr;;1!\Yw;;:f,%
  !\Ygr;;1!\Yw;;;:g,%
  :g:f:d:b:a)
  \draw[ultra thick](0,-4)--(2,-4)--(2,-3)--(3,-3)--(3,-2)--(5,-2)--(5,1);
  \node at (2.5,-5.5) {$ $}; 
\end{tikzpicture}\\
  a\,b^2\,d^3\,(da)^2(fd)^3\,g^2(gf) &
  b\,da\,fd\,gf &
  a^2\,ba\,da\,fd\,gf
\end{array}
$}
\caption{The left diagram is a Chinese staircase with the associated Dyck
  path, and its reading below; the middle its $2$-Chinese staircase and reading; and the
  right its $2$-Chinese function, its reading and its equivalent staircase.}
\label{fig-chinese-2-dyck-reading}
\end{figure}
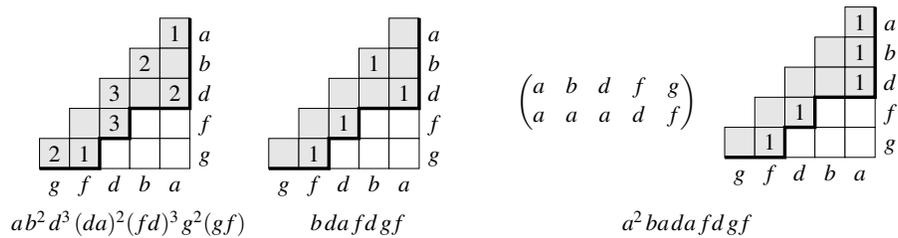

As for the plactic monoid, the elements of this monoid can be represented using a combinatorial object,
the \emph{Chinese staircase}, and its insertion algorithm. It has a particular \emph{row reading}
(defined in~\cite{CEKNH01})
which is the shortlex normal form of its class.
For a chinese staircase $S$,
one can also associate a Dyck path $\dyck(S)$ of length $2|cont(S)|$;
see~\cref{fig-chinese-2-dyck-reading}.

The following result gave us a nice description of the $\chinese(\A,2)$-equivalence.

\begin{thm}\label{thm-chinese-to-2-chinese}
    If $S$ and $T$ are Chinese staircases, then $\RW(S) \equiv_{\chinese(\A, 2)}
    \RW(T)$ if and only if $\cont(\RW(S)) = \cont(\RW(T))$ and $\dyck(S) = \dyck(T)$.
\end{thm}

Using this, we define two different combinatorial objects,
the 2-Chinese staircases and the 2-Chinese functions, that both represent
$\chinese(\A,2)$-classes, in order to have
a better understanding of $\chinese(\A,2)$.

Let $\B\subseteq\A$ and $D$ be a Dyck path of length $2|\B|$. The \emph{2-Chinese staircase}
associated to $\B$ and $D$ is the Chinese staircase $S$ such that $\dyck(S)=D$, that
has 1s in all of the peaks of $D$ and in every boxes of the diagonal having an empty hook.

If $S$ is a 2-Chinese staircase, then the \emph{row reading} $\RW(S)$ of $S$
is simply the row-reading of $S$ as a Chinese staircase which is also the
shortlex normal form of its $\chinese(\A,2)$-class;
see~\cref{fig-chinese-2-dyck-reading}.

A \emph{$2$-Chinese function} on $\A$ is a function $\f:\B\to \B$ for some
$\B\subseteq\A$ which for all $x, y\in \B$ satisfies: $x\leq y$ implies
$\f(x)\leq \f(y)$; and  $\f(x)\leq x$.
We denote by $\emptyf$ the unique function whose domain is empty.

The \emph{insertion} of $y\in\A$ in a $2$-Chinese function $\f$ is
the function $\f\leftarrow y$ whose domain is $\dom(\f)\cup\{y\}$ and
the image of any $x$ is given by $(\f\leftarrow y)(x) := \min\{\{y\} \cup
\f(\widehat{x})\}$, for $\widehat{x}=\min(z\in\dom(\f)|z\geq x)$.
It is routine to verify that $\f \leftarrow y$ is also a $2$-Chinese function,
as seen in the following example:
\[
  \begin{psmallmatrix}
    a&b&c&e&f\\
    a&b&c&c&f\\
  \end{psmallmatrix}
  \leftarrow g =
  \begin{psmallmatrix}
    a&b&c&e&f&g\\
    a&b&c&c&f&g
  \end{psmallmatrix}
  \quad\text{and}\quad
  \begin{psmallmatrix}
    a&b&c&e&f\\
    a&b&c&c&f\\
  \end{psmallmatrix}
  \leftarrow d =
  \begin{psmallmatrix}
    a&b&c&d&e&f\\
    a&b&c&c&c&d\\
  \end{psmallmatrix}
\]
We define the reading word of $\f$ to be
$\RW(\f):=a_1\,\f(a_1)\cdots a_n\,\f(a_n)$,
where $\dom(\f)=\{a_1<a_2<\cdots<a_n\}$.

Using properties of both combinatorial representatives of $\chinese(\A,2)$,
we proved that the Chinese monoid is of type 1.

These objects also allowed us prove that $\chinese(\A,2)$ is $\J$-trivial.
Unlike the stylic monoid, its $\J$-order is surprisingly not graded.

Similar to the plactic case, we define a set of normal forms the following way:
for any $(\f, e)\in \chinese(\A,2)\times \Com(\A, \sigma)$ such that $\cont(\f)=\cont(e)$,
we define the \emph{row reading $\RW_{\sigma}(\f,e)$ of} $(\f, e)$ by inflating
$\RW(\f)$ putting a suitable exponent on the last occurrence of each letter. 
For example, if $\A = \{a,b,c\}$ and $\sigma$ is constant with value $4$, then
$
\RW_\sigma
$
\scalebox{0.8}{
$\left[
\begin{pmatrix}
a&b&c&e&f\\
a&b&c&c&f\\
\end{pmatrix}, a^2b^3c^2ef^3\right]
$}$= a^2b^3c^2ec^3f^3.
$

\subsection{Cardinality and Idempotents of Monoids of Type 1}
\label{section-isomorphism}

From the definition of a type 1 monoid, one only has to know the combinatorial structure of
the 2-quotient in order to compute the cardinality of the $\sigma$-quotient for any $\sigma$
and to find its idempotents.

\begin{thm}
  Let $\sigma:\A \to \N_{\geq2}$ be arbitrary.
  Then the cardinality of $\M(\A,\sigma)$ is given by
  \begin{equation}
    |\M(\A, \sigma)| = \sum_{\B\subseteq \A} \left(s_{|\B|} \prod_{b\in
        \B}\left(\sigma(b) - 1\right)\right)
  \end{equation}
  where $s_k$ is:
  \begin{enumerate}[label=\rm (\roman*)]
  \item the $k$-th Bell number if $\M=\plax$;
  \item the $k$-th Catalan number if $\M=\chinese$;
  \end{enumerate}
  In particular, $|\M(\A, 2)| = \sum_{k=0}^n \binom{n}{k}s_k$ is the binomial
  transform of the sequence $s_k$ in both cases.
\end{thm}

\begin{prop}
  The monoids $\plax(\A,\sigma)$ and $\chinese(\A,\sigma)$
  contain exactly $2^{|\A|}$ idempotents, one for each
  $\B=\{b_1 < b_2 < \cdots < b_k\}\subseteq\A$.
  These elements are
  inflations $I$ of, respectively,
  $\min_{\J}(\plax(\B,2))$ and $\min_{\J}(\chinese(\B,2))$,
  such that $\ev_{\sigma}(I)=\prod_{i=1}^k b_i^{\sigma(b_i)-1}$.
\end{prop}

\section{Monoids of Type 2}
\label{section-Type-2}

Let $\M(\A)$ be an evaluation-preserving monoid. 
We say that $W\in \A ^ *$ is a \emph{gathered element}
if whenever $W'\in \A ^*$ is such that $W'\equiv_{\M(\A)} W$, then $
W' = U'a^2V'$ implies $W = Ua^ 2V$, for some $U, V\in \A ^*$ where
$U$ and $U'$ have the same number of $a$s.

If an $\equiv_{\M(\A)}$-class contains a gathered element, we refer to the
lexicographically largest gathered word $G(W)$ in the class as the
\emph{canonical gathered element}.

We define the \emph{$(a,i)$-expansion} of $W\in A ^ *$ to be the word
obtained from $W$ where the $i$-th occurrence of $a$ is duplicated.

\begin{defi}\label{def:type2}
  An evaluation-preserving monoid $\M(\A)$ is of \emph{type 2} if:
  \begin{enumerate}[label=(\alph*)]
    \item each $\equiv_{\M(\A)}$-class has a gathered element; and
    \item the $(a,i)$-expansion of $G(W)$ equals the canonical gathered element
      of the $(a,i)$-expansion of $W$.
  \end{enumerate}
\end{defi}

If $W\in\A^*$, then we define $G_\sigma(W)$, the \emph{$\sigma$-reduced word of} $W$, to be the word obtained from $G(W)$
by repeatedly replacing any factor $a^{\sigma(a)}$ of $G(W)$ by $a$, until
there are no such factors remaining. The set of $\sigma$-reduced word constitutes a set
of normal forms for $\M(\A,\sigma)$.

\subsection{The $\sigma$-Sylvester Monoids}
\label{section-sylvestre}

The \emph{sylvester
monoid}~\cite{HNT05} is the quotient of $\A^*$ by the following infinite set of
relations:
\[
	acWb = caWb \quad\text{if } a \leq b < c,\quad
	\text{ for all } W\in \A^*.
\]

A \emph{binary search tree} $T = (L, r, R)$ is a binary tree labelled by $\A$
such that the label of each node is greater than or equal to all labels in its
left subtree $L$ and strictly smaller than all labels in its right subtree
$R$.
Binary search trees are endowed with a well-known left insertion $a\rightarrow T$
of letters $a\in\A$;
see \cite{HNT05}.
Given a binary search tree $T$, denote $\RW(T)$ its \emph{right to left postfix reading}.
This word is the
lexicographic largest word in its $\equiv_{\Sylv(\A)}$-class~\cite[Proposition 14]{HNT05}.

\begin{figure}[t]
\centering
\scalebox{0.8}{
\begin{forest}
for tree={grow=south, minimum size=2ex, inner sep=0, l=0}
[,no edge[$2\rightarrow$,no edge[,no edge[,no edge[,no edge]]]]]
\end{forest}
\begin{forest}
for tree={grow=south, circle, draw, minimum size=3ex, inner sep=0, l=0}
[2[2
  [1[,no edge, draw=none]
  	[,no edge, draw=none]]
    [,no edge, draw=none]]
  [3]]
\end{forest}
\begin{forest}
for tree={grow=south, minimum size=2ex, inner sep=0, l=0}
[,no edge[\text{=},no edge[,no edge[,no edge[,no edge]]]]]
\end{forest}
\begin{forest}
for tree={circle, draw, minimum size=3ex, inner sep=0,  l=0}
[2[2
  [1[,no edge, draw=none]
  	[2]]
    [,no edge, draw=none]]
  [3]]
\end{forest}
\hspace{2cm}

\begin{forest}
for tree={grow=south, minimum size=2ex, inner sep=0, l=0}
[,no edge[$2\rightarrow$,no edge[,no edge[,no edge]]]]
\end{forest}
\begin{forest}
for tree={grow=south, circle, draw, minimum size=3ex, inner sep=0, l=0}
[2
  [1[,no edge, draw=none]
  	[,no edge, draw=none]]
  [3]]
\end{forest}
\begin{forest}
for tree={grow=south, minimum size=2ex, inner sep=0, l=0}
[,no edge[\text{=},no edge[,no edge[,no edge]]]]
\end{forest}
\begin{forest}
for tree={grow=south, circle, draw, minimum size=3ex, inner sep=0, l=0}
[2[1[,no edge, draw=none]
  	[2]]
  [3]]
\end{forest}
}
\caption{\label{insere-3}On the left, the insertion of $2$ in the binary search tree that has reading word $3122$ resulting to the tree having reading word $32122$. On the right, the 2-sylvester insertion in its 2-reduced binary search tree.}
\vspace{-.7em}
\end{figure}
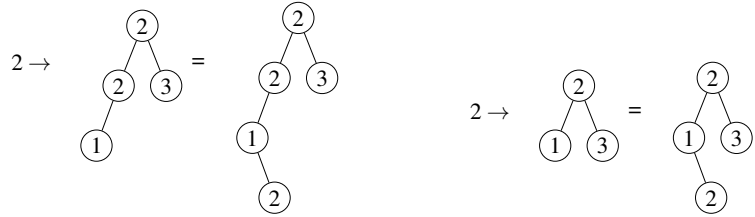

Using this set of normal forms, one can prove that the sylvester monoid is of type 2.

The $2$-reduced words are the
readings of the binary trees where parents have different label than their left children.
If the number of nodes $i$, and the number of nodes $k$
having at least one ancestor with the same label are fixed, then,
thanks to~\cite{CY19} and a standard involution among trees,
\emph{recursive reversal of the left branch subtrees},
one can prove that the number of such trees is $B_{i-k-1,k}$, an element of the Borel
triangle (A234950 in~\cite{oeis}).

\begin{thm}\label{EnumBST}
	Let $\A$ be an alphabet of size $n\in \N$. Then
	\begin{equation}
    |\Sylv(\A, 2)| = 1 +
    \sum_{i=1}^{2n-1} \sum_{k=0}^{\lfloor i/ 2\rfloor}
    B_{i-k-1,k} \binom{n}{i-k}.
  \end{equation}
\end{thm}

One can easily adapt the enumeration formula for $\Sylv(\A,p)$ but,
not being of type 1 makes the general formula 
significantly more complicated.

\begin{prop}
The monoid $\Sylv(\A,2)$ contains exactly $\sum_{k=0}^n \binom{n}{k} S_k$ idempotents,
where $S_n$ is the  Schr\"oder numbers
(A006318 in~\cite{oeis}), shifted by one:
$S_0=S_1=1$, $S_2=2$, $S_3=6$, \emph{etc}.

These idempotents are the
postfix reading of
2-reduced binary search trees $T$ such that,
for all $x\in\cont(T)$, the deepest node labelled $x$ in $T$ does not have any left subtree;
see the rightmost tree of \cref{insere-3}.
\end{prop}

Using these trees, one can describe the idempotents of $\Sylv(\A,\sigma)$ for arbitrary $\sigma$
but the enumeration formula can only be easily adapted for $\sigma$ constant.
\vspace{-1.6em}

\bibliographystyle{eptcs}
\bibliography{main}

\end{document}